\documentclass[11pt,twoside,a4paper]{amsart}
\usepackage{amsmath,amssymb,amsthm}

\usepackage{graphicx,psfrag}
\usepackage{amssymb}

\newtheorem{Thm}{Theorem}[section]
\newtheorem{Lem}[Thm]{Lemma}
\newtheorem{Pro}[Thm]{Proposition}
\newtheorem{Cor}[Thm]{Corollary}

\theoremstyle{definition}

\theoremstyle{remark}
\newtheorem{Rem}[Thm]{Remark}

\newcommand{\R}{\mathbb{R}}

\newcommand{\ga}{\gamma}

\renewcommand{\phi}{\varphi}

\date{}

\begin{document}

\title{Non-Positive Curvature and the Ptolemy Inequality}
\author{Thomas Foertsch,  
\ Alexander Lytchak, \  Viktor Schroeder$^1$
}

\footnote{$^1$Supported by Swiss National Science Foundation}

\address{Thomas Foertsch, Alexander Lytchak\\ Mathematisches Institut\\ Universit\"at Bonn\\
Beringstr. 1, 53115 Bonn, Germany}
\address{Viktor Schroeder\\ Mathematisches Institut\\ Universit\"at Z\"urich\\
Winterthurerstrasse 190\\ 8057 Z\"urich\\ Switzerland }

\email{foertsch\@@math.uni-bonn.de \\ lytchak\@@math.uni-bonn.de \\ vschroed\@@math.unizh.ch}




\begin{abstract}
We provide examples of non-locally compact geodesic Pto\-lemy metric spaces which are not
uniquely geodesic. On the other hand, we show that locally compact, geodesic Ptolemy metric 
spaces are uniquely geodesic. Moreover, we prove that a metric space is $\operatorname{CAT}(0)$ if and only if
it is Busemann convex and Ptolemy.
\end{abstract}

\maketitle


\section{Introduction}

A metric space $X$ is called a Ptolemy metric space, if the inequality
\begin{equation} \label{eqn-ptolemy}
|xy||uv| \; \le \; |xu| \, |yv| \; + \; |xv| \, |yu|
\end{equation}
is satisfied for all $x,y,u,v\in X$. \\

Our interest in Ptolemy metric spaces 
originates from an analysis of boundaries of 
$\operatorname{CAT}(-1)$-spaces when endowed with a Bourdon or Hamenst\"adt metric.
Such boundaries are indeed Ptolemy metric spaces (see \cite{foesch}). 

Various aspects of such spaces have occasionally been studied, for instance in \cite{bufw}, \cite{k}, \cite{klm}
and \cite{sm}. A smooth Riemannian manifold is of non-positive sectional
curvature, if and only if it is locally 
Ptolemy, and a locally Ptolemy Finsler manifold is necessarily Riemannian (see \cite{k2} and \cite{bufw}).
Furthermore, $\operatorname{CAT}(0)$-spaces are geodesic Ptolemy metric spaces (compare Section \ref{sec-prel}). 

On first consideration, these observations might suggest that for geodesic metric spaces the Ptolemy condition
is some kind of non-positive curvature condition.
We show that without any further conditions this is completely wrong.
\begin{Thm} \label{theo-isom-emb}
Let $X$ be an arbitrary Ptolemy space,
then $X$ can be isometrically embedded into
a complete geodesic Ptolemy space $\hat X$.
\end{Thm}

As an application, take the four
point Ptolemy space
$X=\{x,y,m_1,m_2\}$ with
$|xy|=2$ and all other nontrivial distances
equal to $1$.
By Theorem \ref{theo-isom-emb} $X$ can be isometrically embedded into
a geodesic Ptolemy space $\hat X$. Since 
$m_1$ and $m_2$ are midpoints of $x$ and $y$,
there are in $\hat X$ two different geodesics
joining the points $x$ and $y$.
In particular, the space is 
not uniquely geodesic and hence far away from 
'non-positively curved'. 

The space $\hat X$ constructed in Theorem \ref{theo-isom-emb} fails to be proper. 
Indeed, for proper geodesic Ptolemy spaces 
the situation is completely different:

\begin{Thm} \label{theo-ptolemy-distance-convex}
A proper, geodesic, Ptolemy metric space is uniquely geodesic.
\end{Thm}

The proof of Theorem \ref{theo-ptolemy-distance-convex} also works
if one replaces 'properness' through the assumption that there exists a geodesic
bicombing which varies continuously with its endpoints. This will
help us to prove that the property of being Ptolemy precisely distinguishes between
the two most common non-positive curvature conditions for
geodesic metric spaces, namely the one due to Alexandrov and the one due to Busemann.

\begin{Thm} \label{theo-cat0-busemann-ptolemy}
A metric space is $\operatorname{CAT}(0)$ if and only if it is Ptolemy and Busemann convex.
\end{Thm}

We finish the Introduction with a few comments on the theorems above. \\
The space $\hat{X}$ constructed in Theorem \ref{theo-isom-emb} has several remarkable
properties. There is a huge collection of convex functions on $\hat{X}$, in particular, all
distance functions to points are convex. However, from the geometrical or topological point
of view the space appears rather odd. Starting with finite spaces, for instance
with the four point space $X$ as above, one obtains a space 
$\hat{X}$ on which the distance functions $d_x$ to the points $x\in X\subset \hat{X}$ are affine.
This is to our knowledge a first example of very strange affine functions on metric
spaces. The space $\hat{X}$ is far from being a product, in contrast to the structural
results obtained in \cite{lsch} and \cite{hl}. 
Theorem \ref{theo-isom-emb} shows that the Ptolemy condition is sufficient to ensure that
a metric space can be isometrically embedded into a geodesic Ptolemy metric space. This seems to be
particularly interesting, as the problem of synthetic descriptions of (non-convex) subsets of
$\operatorname{CAT}(0)$-spaces is very difficult (cf. \cite{g}).
The properness, assumed in addition
in Theorem \ref{theo-ptolemy-distance-convex}, forces the space to be contractible; in fact an absolut 
neighborhood retract. Moreover, it can be shown that the distance functions
to points are neither affine in this case, in fact they are strictly convex 
(cf. Remark \ref{rem-gen-strict-dist-convex}). Finally, the proof of 
Theorem \ref{theo-cat0-busemann-ptolemy} confirms the following idea. If a Busemann convex
space is not $\operatorname{CAT}(0)$, then it contains an infinitesimal portion of a non-Euclidean
Banach space. This observation may be of some interest in its own right. \\

After a preliminary section,
we prove Theorem \ref{theo-isom-emb} in Section \ref{sec-theo-isom-emb},
Theorem \ref{theo-ptolemy-distance-convex} in Section 
\ref{sec-theo-ptolemy-distance-convex} and Theorem \ref{theo-cat0-busemann-ptolemy} in
Section \ref{sec-theo-cat0-busemann-ptolemy}.


\section{Preliminaries}
\label{sec-prel}

We start by recalling some easy examples.
The real line $\R$ is a Ptolemy space. To show this,
consider points
$x,y,z,w$ in this order on the line.
A completely trivial computation shows that $|xz| |yw| =|xy||zw| + |yz| |wx|,$
which implies that $\R$ is Ptolemy.
Since the Ptolemy condition on four points 
is M\"obius invariant (see below),
the equality 
above holds for points 
$x,y,z,w$, which lie in this order on a circle in the
plane $\R^2$. This is the classical Theorem of
Ptolemy for cyclic quadrilaterals.

To show that the Euclidean space $\R^n$ is Ptolemy,
consider again four points
$x,y,z,w$. Applying a suitable M\"obius transformation
we can assume that
$z$ is a midpoint of $y$ and $w$, i.e.
$|yz|=|zw|=\frac{1}{2}|yw|$.
For this configuration the
Ptolemy inequality is equivalent
to 
$|xz| \le \frac{1}{2}(|xy|+|xw|)$,
which is just the convexity of the distance
to the point $x$.

Every $\operatorname{CAT}(0)$-space is Ptolemy, since every four point
configuration in a $\operatorname{CAT}(0)$-space admits a 
subembedding into the Euclidean plane
(\cite{brh}, p. 164). \\



{\bf M\"obius invariance:}
Let $d, d'$ be two metrics on the same set
$X$. The metrics are called M\"obius equivalent,
if for all quadruples $x,y,z,w$ of points
$$\frac{d(x,y)d(z,w)}{d(x,z)d(y,w)} =
\frac{d'(x,y)d'(z,w)}{d'(x,z)d'(y,w)} .$$

If $d$ and $d'$ are M\"obius equivalent, then
$(X,d)$ is Ptolemy, if and only if
$(X,d')$ is Ptolemy.

Indeed, the Ptolemy inequality says that
for all quadruples of points the three numbers  of the triple
$$A=(d(x,y)d(z,w),d(x,z)d(y,w),d(x,w)d(y,z))$$
satisfy the triangle inequality.
By dividing all three numbers by $d(x,w)d(y,z)$,
we see that these numbers satisfy the triangle inequality if
and only if the three numbers of the triple
$$B=(\frac{d(x,y)d(z,w)}{d(x,w)d(y,z)},\frac{d(x,z)d(y,w)}{d(x,w)d(y,z)},1)$$
satisfy the triangle inequality. In this expression we can replace
$d$ by $d'$, and hence we obtain the claim. \\



{\bf Basic properties of Ptolemy spaces:}
Here we collect a couple of basic properties of Ptolemy spaces which will be frequently used in the
remainder of this paper.
\begin{description}
\item[(P1)] Every subset $Y\subset X$ of a Ptolemy metric space 
$X$, endowed with the metric inherited from $X$, is Ptolemy.
\item[(P2)] A metric space $X$ is Ptolemy if and only
if for every $\lambda > 0$ the scaled space $\lambda X$ is Ptolemy.
\end{description}

Some of our arguments below will use the notions of ultrafilters and ultralimits;
a generalization of pointed Gromov-Hausdorff convergence. 
We refer the reader, not familiar with these methods, to \cite{brh} and \cite{l}. 
The symbol $\lim_\omega (X_n,x_n)$ will denote such an ultralimit (w.r.t. a non-principal ultrafilter $\omega$). 

As every metric property, the Ptolemy condition is invariant w.r.t. ultraconvergence.
\begin{description}
\item[(P3)] For every sequence $\{ (X_i,x_i)\}_i$ of pointed Ptolemy spaces and every non-principle 
ultrafilter $\omega$, the ultralimit $\lim_\omega (X_i,x_i)$ is a Ptolemy space.
\end{description}
Finally, we recall another important observation, which is due to Schoenberg (see \cite{sch}).
This property lies in the heart of Theorem \ref{theo-cat0-busemann-ptolemy}.
\begin{description}
\item[(P4)] A normed vector space is an inner product space if and 
only if it is Ptolemy. 
\end{description} 

A subset of a normed vector space is called linearly convex, if with any two points
it contains the straight line segment connecting these points. A metric space is called
linearly convex, if it is isometric to a linearly convex subset of a normed vector space and
called flat, if it is isometric to a convex subset of an inner product space. \\
With this notation the properties above immediately yield the

\begin{Cor} \label{cor:schoenberg}
Let $X$ be a Ptolemy space, then every linearly convex subset $C\subset X$ of $X$ is flat. \\
\end{Cor}



Let $X$ be a metric space and $x,y\in X$. Then a point $m\in X$ is called a midpoint of $x$ and $y$
if $|xm|=\frac{1}{2}|xy|=|my|$. We say that $X$ has the midpoint property, if for all $x,y\in X$ there exists
a midpoint of $x$ and $y$ in $X$.
The space $X$ is called geodesic, if any two points $x,y\in X$ can be joined by a geodesic path, i.e. 
to any two points $x,y\in X$ there exists an isometric embedding $\gamma$ of the interval $[0,|xy|]$
of the Euclidean line into $X$ such hat $\gamma (0)=x$ and $\gamma (|xy|)=y$.
In this paper we will always assume that geodesics are parameterized
affinely, i.e. proportionally to arclength. A complete space with the midpoint property is geodesic.
Note that a Ptolemy metric space $X$ which has the midpoint property satisfies
\begin{equation} \label{eqn-def-distance-convexity}
|mz| \; \le \; \frac{1}{2} \, [|xz| \, + \, |yz|] 
\end{equation}
for all $x,y,z,m\in X$ such that $m$ is a midpoint of $x$ and $y$ (cf. \cite{foe} for a 
discussion of such spaces). \\



{\bf Non-positive curvature conditions and often convex spaces:}
The most common non-positive curvature conditions are due to Alexandrov and Busemann 
(cf. \cite{a} and \cite{bus}). \\
We suppose that the reader interested in this paper's subject is familiar with the notion of 
$\operatorname{CAT}(0)$-spaces.
Roughly speaking, a geodesic space $X$ is called a $\operatorname{CAT}(0)$-space
if all geodesic triangles in $X$ are not thicker than their comparison triangles in the
Euclidean plane $\mathbb{E}^2$ (for a precise definition we refer the reader, for instance, to
Section II.1 in \cite{brh}). \\
The geodesic space $X$ is said to be Busemann convex, if 
for any two (affinely parameterized) geodesics $\alpha ,\beta :I\longrightarrow X$, the map
$t\mapsto |\alpha (t)\beta (t)|$ is convex. 

\begin{Rem}
Every $\operatorname{CAT}(0)$-space is Busemann convex and every Busemann convex space is
uniquely geodesic. However, there are Busemann convex
spaces which are not $\operatorname{CAT}(0)$-spaces, as for instance all non-Euclidean normed vector spaces 
with strictly convex unit norm balls. 
\end{Rem}

The notion of Busemann convexity is not stable under limit operations. For instance,
a sequence of strictly convex norms on $\R^n$ may converge to
a non-strictly convex norm.
To overcome such a
phenomenon, Bruce Kleiner introduced a weaker notion of 
Busemann convexity in \cite{kl},
which is stable with respect to limit operations. 
A metric space $X$ is called 
{\em often convex}, if there exists
a convex geodesic bicombing, i.e. a map
$\ga:X\times X \times [0,1] \to X$,
$(x,y,t)\mapsto \ga_{x,y}(t)$, such that
$t\mapsto  \ga_{x,y}(t)$ is a geodesic and for all
$x,y,x',y'$ and all
(not necessarily surjective) affine maps
$\phi:I \to [0,1]$,
$\psi:I \to [0,1]$ defined on the same interval $I$ the
map
$s \mapsto |\ga_{x,y}(\phi(s))\ga_{x'y'}(\psi(s))|$
is convex. 
This convexity implies in particular, that
$\ga$ is continuous and that
for
points
$x',y' \in \ga_{x,y}([0,1])$ the
geodesic
$\ga_{x',y'}$ is contained in $\ga_{x,y}$.
In particular, $\ga$ defines a continuous midpoint map
$m(x,y):=\ga_{x,y}(1/2)$.

We state the following properties of
often convex spaces, that are direct consequences of the definitions (cf. \cite{kl}).

\begin{description}
\item[(OC1)] 
For every sequence $\{ (X_i,x_i)\}_i$ of often convex spaces and every 
non-principle ultrafilter $\omega$, 
the ultralimit $\lim_\omega (X_i,x_i)$ is an often convex space.
\item[(OC2)] A metric space $X$ is often convex if and only
if for every $\lambda > 0$ the scaled space $\lambda X$ is often convex.
\item[(OC3)] Let $X$ be often convex. 
Then $X$ is Busemann convex if and only if $X$ is uniquely geodesic.
\end{description}


\section{Proof of Theorem \ref{theo-isom-emb}}
\label{sec-theo-isom-emb}

We explicitly construct the complete, geodesic, Ptolemy metric space $\hat X$.
First, we subsequently add midpoints to $X$ in order to obtain a Ptolemy
metric space $\mathfrak{M}(X)$ which has the midpoint property. Then we pass to an ultraproduct
of $\mathfrak{M}(X)$. \\
Let $\Sigma$ denote the set of 
unordered tuples in $X$. 
Formally,
$\Sigma= \{ \{x_1,x_2\} \subset X| x_1,x_2\in X\}$,
i.e. $\Sigma$ consists of all subsets of $X$ with one or two elements.

On $\Sigma$ we define a metric via
\begin{displaymath}
|\{x_1,x_2\}\{y_1,y_2\}| \; := \; 
\left\{
\begin{array}{cl}
\frac{1}{4} [|x_1y_1|+|x_1y_2|+|x_2y_1|+|x_2y_2|] & \mbox{if} \; \{x_1,y_1\} \neq \{y_1,y_2\} \\
0 & \mbox{if} \; \{x_1,x_2\}=\{y_1,y_2\}
\end{array}
\right\}
\end{displaymath}
for all $\{x_1,x_2\},\{y_1,y_2\}\in \Sigma$. This indeed defines a metric on $\Sigma$. In order to verify this,
one has to prove the triangle inequality
\begin{displaymath}
|\{ x_1,x_2\}\{ y_1,y_2\}| \; \le \; |\{ x_1,x_2\} \{ z_1,z_2\} | \; + \; |\{ z_1,z_2\} \{ y_1,y_2\} | 
\hspace{0.5cm} 
\end{displaymath}
for all $\{ x_1,x_2\} ,\{ y_1,y_2\} ,\{ z_1,z_2 \}\in \Sigma$.
If two of the tuples coincide, the validity of the inequality above is evident, and
otherwise it just follows by applying the triangle inequality in $X$ three times. \\
Moreover, the space $M(X):=(\Sigma ,|\cdot|)$ is Ptolemy, i.e., it satisfies
\begin{align*}
|\{ x_1,x_2\} \{ y_1,y_2\} | \cdot |\{ z_1,z_2\} \{ u_1,u_2\} | &
\le |\{x_1,x_2\} \{ z_1,z_2\} | \cdot |\{ y_1,y_2\} \{ u_1,u_2\} | \\ &
+  |\{x_1,x_2\} \{ u_1,u_2\}| \cdot |\{ y_1,y_2\} \{ z_1,z_2\} |
\end{align*}
for all $\{ x_1,x_2\},\{ y_1,y_2\} ,\{ z_1,z_2\} ,\{ u_1,u_2\} \in \Sigma$. Once again, the 
validity of this inequality is evident, if
two of the tuples coincide, and otherwise it follows by applying the Ptolemy inequality 
in $X$ sixteen times. \\
Note further that $X$ isometrically embeds into $M(X)$ via $x\mapsto \{x,x\}$. 
Thus we may identify $X$ with
a subset of $M(X)$. \\
Now we define $M^0(X):=X$ as well as $M^{k+1}(X):=M(M^k(X))$ and set 
$\mathfrak{M}(X):=\bigcup\limits_{k=0}^\infty \, M^k(X)$. From the considerations above it follows that this space 
is a Ptolemy metric space. Moreover, it has the midpoint property. 
Namely, each pair $x,y\in \mathfrak{M}(X)$ is contained in some $M^k(X)$ and $\{ x,y\} \in M^{k+1}(X)$
is a midpoint of $x$ and $y$.
Passing to an ultraproduct $\hat{X}$ of 
$\mathfrak{M}(X)$, i.e. $\hat{X} := \lim_\omega \{ (\mathfrak{M}(X), \mathfrak{x})\}_n$, the ultralimit of 
the constant sequence 
$\{ (\mathfrak{M}(X), \mathfrak{x})\}_n$ w.r.t. some ultrafilter $\omega$, where $\mathfrak{x}\in \mathfrak{M}(X)$, 
we obtain a complete Ptolemy metric space which has the midpoint property and, therefore, is geodesic.
\hfill $\Box$


\section{Proof of Theorem \ref{theo-ptolemy-distance-convex}}
\label{sec-theo-ptolemy-distance-convex}

In this section we prove Theorem \ref{theo-ptolemy-distance-convex}. For two points $p^-,p^+\in X$
we consider the set
\begin{displaymath}
C(p^-,p^+) \; := \; \Big\{ x\in X \, \Big| \, |p^-p^+| \, = \, |p^-x| \, + \, |xp^+| \Big\} .
\end{displaymath}
Since the distance functions 
$d_{p^\pm}:=d(p^\pm ,\cdot):X \rightarrow \mathbb{R}_0^+$ are convex,
the set $C(p^-,p^+)$ is convex.
The convexity of
$d_{p^\pm}$ implies that these functions are
affine on the set
$C(p^-,p^+)$. \\

{\bf Proof of Theorem \ref{theo-ptolemy-distance-convex}:}
Let $\gamma_1,\gamma_2:[0,L]:=|p^-p^+|\rightarrow C$ be  
geodesics connecting $p^-$ to $p^+$. Set $x_s:=\gamma_1(s)$ as well as $y_s:=\gamma_2(s)$ and let
$m_s$ denote a midpoint of $x_s$ and $y_s$. Since the 
functions $d_{p^\pm}:C(p^-,p^+)\rightarrow \mathbb{R}_0^+$ are affine, we obtain
$|p^-m_s|=s$ and $|p^+m_s|=L-s$. \\
Let now $0<s<t\le L$ be arbitrary, then 
the triangle inequality yields
\begin{equation}\label{eq:sec4.1}
|m_sx_t|+|m_sy_t|\ge |x_ty_t|,
\end{equation}
the Ptolemy inequality yields
\begin{displaymath}
|p^-m_t| \cdot |x_sy_s| \; \le \; |m_tx_s|\cdot |p^-y_s| \; + \; |m_ty_s|\cdot |p^-x_s|
\end{displaymath}
and therefore
\begin{equation}\label{eq:sec4.2}
|m_tx_s| \; + \; |m_ty_s| \; \ge \; \frac{t}{s} \cdot |x_sy_s| \; = \; |x_sy_s| \: + \; (t-s) \frac{|x_sy_s|}{s} .
\end{equation}
On the other hand the Ptolemy property
for the points $m_s,x_t,m_t,x_s$ and $m_s,y_t,m_t,y_s$ yields
\begin{equation}\label{eq:sec4.3}
|m_sx_t|\cdot |m_tx_s| \; + \; |m_sy_t|\cdot |m_ty_s| \; \le \;
\frac{1}{2} |x_sy_s|\cdot |x_ty_t| \; + \; 2 (t-s) \cdot |m_sm_t|.
\end{equation}
Now, fix any $0<s< L$ with $l:=|x_sy_s|\ge 0$, and choose a sequence $s_n\rightarrow s$ with $s_n>s$. Let
$m_s$ denote midpoints of $x_{s_n}$ and $y_{s_n}$. By compactness, we can pass to a subsequence 
$m_{s_n}\rightarrow m$, where $m=m_s$ is a midpoint of $x_s$ and $y_s$. Now, we set
\begin{displaymath}
\begin{array}{lllll}
\varphi_n := |m_{s_n}m|, & \epsilon_n := s_n-s, & a_n^+:=|m_{s_n}x_s|, & b_n^+:=|m_{s_n}y_s| & \\
a_n^-:= |m,x_{s_n}|, & b_n^- := |my_{s_n}|, & l_n := |x_{s_n}y_{s_n}| & \mbox{and} & Q:=\frac{l}{s}\ge 0.
\end{array}
\end{displaymath} 

\begin{figure}[htbp]
\psfrag{m}{$m$}
\psfrag {msn}{$m_{s_n}$}
\psfrag{l}{$l$}
\psfrag{ln}{$l_n$}
\psfrag{epsn}{$\epsilon_n$}
\psfrag{xs}{$x_s$}
\psfrag{ys}{$y_s$}
\psfrag{xsn}{$x_{s_n}$}
\psfrag{ysn}{$y_{s_n}$}
\psfrag{a+n}{${\scriptstyle a^+_n}$}
\psfrag{a-n}{${\scriptstyle a^-_n}$}
\psfrag{b+n}{${\scriptstyle b^+_n}$}
\psfrag{b-n}{${\scriptstyle b^-_n}$}
\psfrag{p-}{$p^-$}
\psfrag{p+}{$p^+$}
\psfrag{phin}{${\scriptstyle \varphi_n}$}
\includegraphics[width=0.9\columnwidth]{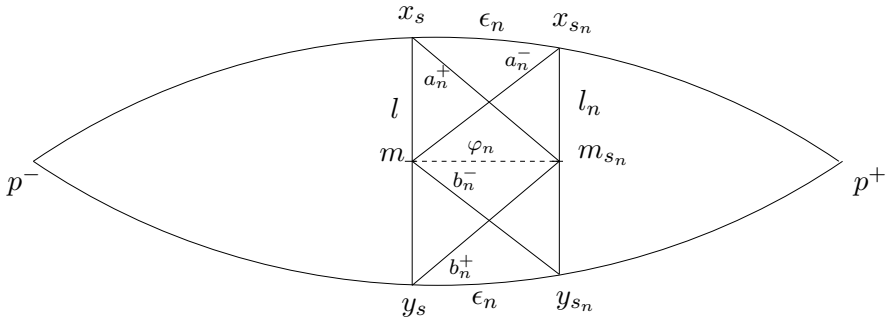}
\caption{This figure visualizes the notation used in the proof.}
\end{figure}

We have by (\ref{eq:sec4.1}), (\ref{eq:sec4.2}) and 
 (\ref{eq:sec4.3}):
\begin{equation}\label{eq:sec4.4}
a_n^- +b_n^-\ge l_n, \;\;\;
a_n^+ +b_n^+\ge l+\epsilon_n \cdot Q, \;\;\;
b_n^+ \cdot b_n^- + a_n^+\cdot a_n^- \le \frac{1}{2}l\cdot l_n+2\epsilon_n\varphi_n .
\end{equation}

By triangle inequalities, we see
$|l_n-l|\le 2\epsilon_n$, $|a^\pm_n-\frac{1}{2}l|\le 2\epsilon_n$,
$|b^\pm_n-\frac{1}{2}l|\le 2\epsilon_n$, hence
passing to a subsequence, we find $A^+,A^-,B^+,B^-,C\in [-2,2]$ such that
$a^\pm_n =\frac{1}{2}l+A^\pm \epsilon_n+o(\epsilon_n )$, $b^\pm_n =\frac{1}{2}l+B^\pm \epsilon_n+o(\epsilon_n )$ and
$l_n=l+C\epsilon_n+o(\epsilon_n)$. 
Since $2\epsilon_n\varphi_n = o(\epsilon_n)$, we
derive from (\ref{eq:sec4.4}) at
\begin{displaymath}
B^-+A^-\ge C,\;\; B^++A^+\ge Q\ge 0,\;\; 
(B^++B^-)+(A^++A^-) \le C. 
\end{displaymath}
Thus $Q = 0$, which implies
$\ga_1(s)=\ga_2(s)$.
Since $s\in (0,L)$ was arbitrary, we have $\ga_1=\ga_2$. \hfill $\Box$

\begin{Rem} \label{rem-gen-strict-dist-convex}
It is not difficult to prove that in Theorem \ref{theo-ptolemy-distance-convex}
one can replace ''properness'' by the slightly weaker assumption of ''local compactness''. 
A more involved argument shows that a locally compact geodesic Ptolemy metric space is even
strictly distance convex, i.e. that the inequality (\ref{eqn-def-distance-convexity})
is strict, whenever $|xy|>||xz|-|zy||$. The proof of this claim will be given elsewhere.
\end{Rem}

We finish this section with an important observation. In the proof 
of Theorem \ref{theo-ptolemy-distance-convex} we use the properness of $X$
only to show that
$m_{s_n}\to m$.
Clearly, the existence of a continuous midpoint map
$m:X \times X \to X$ implies also this convergence.
As a consequence we have

\begin{Thm} \label{theo-cont-bicombing}
Let $X$ be a geodesic, Ptolemy space which admits a 
continuous midpoint map. Then 
$X$ is uniquely geodesic.
\end{Thm}

Since often convex spaces admit continuous midpoint maps 
the first statement of the following corollary is 
an immediate consequence.

\begin{Cor} \label{cor-often-ptolemy:busemann}
Let $X$ be an often convex and Ptolemy space.
Then $X$ is uniquely geodesic and hence, by (OC3),
Busemann convex.
Moreover,
for every sequence of pointed Busemann convex and 
Ptolemy metric spaces $\{ (X_i,x_i)\}$, every non-principle ultrafilter
$\omega$ and all $\lambda_i \in \mathbb{R}^+$, 
the ultralimit $\lim_{\omega} (\lambda_i X_i,x_i)$
is Busemann convex. 
\end{Cor}

The second statement goes as follows:
A Busemann convex space is in particular often convex.
By (OC1) this property and by (P3) the Ptolemy property
is stable under limits. Thus the limit is
often convex and Ptolemy and hence Busemann convex.


\section{Proof of Theorem \ref{theo-cat0-busemann-ptolemy}}

\label{sec-theo-cat0-busemann-ptolemy}

In this section we prove Theorem \ref{theo-cat0-busemann-ptolemy}. First, we
recall the notion of generalized and weak angles (cf. \cite{l}). Then 
we prove Proposition \ref{prop-convex-hull}
(a version of the Toponogov rigidity theorem), that was shown independently by Rinow and Bowditch
under slightly more restrictive conditions (see \cite{r} and \cite{bow}), and which will finally allow us to
obtain a proof of our Theorem \ref{theo-cat0-busemann-ptolemy}.



\subsection{Ultrarays associated to geodesics and their enclosed angles}

{\bf Ultrarays:} Let $X$ be a geodesic metric space and let $\gamma$ be a geodesic in $X$
emanating from $p\in X$. Now take a non-principle ultrafilter $\omega$, consider the $\omega$-blow up
$(\bar{X},\bar{d})$ of $X$ in $p$, i.e. $(\bar{X},\bar{d}):=\lim_\omega \{ (nX,p)\}_n$,
and define $\bar{\gamma}:[0,\infty )\longrightarrow \bar{X}$ through
$\bar{\gamma}(s) \; := \; \lim_\omega \{ (\gamma (\frac{s}{n}),p)\}_n$ for all $s\in [0,\infty )$.
This map indeed is a geodesic ray in $(\bar{X},\bar{d})$ emanating in $\{ p\}_n\in \bar{X}$. 
We call $\bar{\gamma}$ the ultraray associated to $\gamma$ (and $\omega$). \\

{\bf Weak angles:} In order to get a grip on the interplay between geodesics and their associated 
ultrarays, we recall certain notions of angles. \\
Given three points $p$, $x$ and $y$ in a metric space $X$, consider corresponding comparison points 
$p'$, $x'$ and $y'$ in the Euclidean plane $\mathbb{E}^2$. Let $[p',x']$ and $[p',y']$ denote the 
geodesic segments in $\mathbb{E}^2$ connecting $p'$ to $x'$ and $p'$ to $y'$. These segments enclose 
an angle in $p'$ and this angle is referred to as the (Euclidean)  comparison angle of $x$ and $y$ at $p$.
We write $\angle_p(x,y)$ for this angle. \\
Let now $X$ be a metric space and consider two geodesic segments $\gamma_1$ and $\gamma_2$
parameterized by arclength,
both initiating in some $p\in X$. Then $\gamma_1$ and $\gamma_2$ are said to enclose the angle
$\angle_p(\gamma_1, \gamma_2)$ (in the strict sense) at $p$ if the limit 
$\angle_p(\gamma_1, \gamma_2):=\lim\limits_{s,t\to 0}  \angle_p(\gamma_1(s),\gamma_2(t))$
exists. \\
Recall that for instance a normed vector space is an inner product space if and only if all straight line
segments emanating from the origin enclose an angle. However, even in normed vector spaces that are not inner
product spaces certain so called generalized angles do exist between any straight line segments initiating 
in a common point. Such generalized angles were introduced in \cite{l}. \\
Let $a,b>0$ and $\gamma_1$ and $\gamma_2$ be as above. then we say that $\gamma_1$ and $\gamma_2$
enclose a generalized angle $\angle_p^g(\gamma_1,a,\gamma_2,b)$ at scale $(a,b)$, if the limit
\begin{displaymath}
\angle_p^g(\gamma_1,a,\gamma_2,b) \; := \; \lim\limits_{s\rightarrow 0} \, \angle_p (\gamma_1(as),\gamma_2(bs))
\end{displaymath}
exists. If $\gamma_1$ and $\gamma_2$ enclose generalized angles at all scales $(a,b)$ and,
moreover, these generalized angles do not depend on the particular scale, then we say that $\gamma_1$
and $\gamma_2$ enclose the weak angle
\begin{displaymath}
\angle_p^w (\gamma_1,\gamma_2) \; := \; \angle_p^g (\gamma_1,1,\gamma_2,1).
\end{displaymath}

{\bf Ultrarays in Busemann convex spaces:} Now let $X$ be Busemann convex and $\gamma_1$ and $\gamma_2$
be geodesics on $X$ with $\gamma_1(0)=p=\gamma_2(0)$. Then, for all scales 
$(a,b)$, the generalized angle $\angle_p^g(\gamma_1,a,\gamma_2,b)$ exists. This is immediate, since
$s\mapsto \frac{|\gamma_1(as)\gamma_2(bs)|}{s}$ is monotonously increasing. \\
Next consider the ultrarays $\overline{\gamma_1}$ and $\overline{\gamma_2}$ associated to $\gamma_1$ and $\gamma_2$.
These ultrarays satisfy
\begin{equation} \label{eqn-ultrarays}
\bar{d} \Big( \overline{\gamma_1}(as),\overline{\gamma_2}(bs)\Big) \; = \; 
s \cdot \bar{d}\Big( \overline{\gamma_1}(a),\overline{\gamma_2}(b) \Big) \hspace{1cm} \forall a,b,s>0.
\end{equation}
Moreover, the existence of weak angles of geodesics $\gamma_1$ and $\gamma_2$ in a Busemann
convex space is equivalent to the existence of angles (in the strict sense) between their
associated ultrarays $\overline{\gamma_1}$ and $\overline{\gamma_2}$ in $\bar{X}$.

\begin{Lem} \label{lemma-rays-ultrarays-angles}
Let $X$ be Busemann convex, let $\gamma_1$ and $\gamma_2$ denote geodesics in $X$ initiating in
a common point $p\in X$ and let  $\overline{\gamma_1}$ and $\overline{\gamma_2}$ denote
their associated ultrarays. Then the following properties are mutually equivalent.
\begin{enumerate}
\item The rays $\gamma_1$ and $\gamma_2$ enclose a weak angle.
\item The ultrarays  $\overline{\gamma_1}$ and $\overline{\gamma_2}$ enclose an angle (in the strict sense).
\item The union  $\overline{\gamma_1}(\mathbb{R}^+) \cup \overline{\gamma_2}(\mathbb{R}^+)$ admits
an isometric embedding into the Euclidean plane $\mathbb{E}^2$.
\end{enumerate} 
\end{Lem}
\begin{proof}
The equivalence of (2) and (3) follows immediately from Equation \ref{eqn-ultrarays}.
Moreover, this equation also implies that the ultrarays $\overline{\gamma_1}$ and $\overline{\gamma_2}$
enclose an angle if and only if they enclose a weak angle. Hence the equivalence of (1) and (2) is
a consequence of 
$\frac{1}{s}\bar{d} ( \overline{\gamma_1}(as),\overline{\gamma_2}(bs))=
\lim\limits_\omega \{\frac{1}{n}|\gamma_1(a/n)\gamma_2(b/n)|\}_n$,
and the fact that the generalized angles between $\gamma_1$ and $\gamma_2$ exist for all scales in any case.
\end{proof}

\begin{Pro} \label{prop-cat0-angles}
Let $X$ be Busemann convex.
Assume that for all geodesic segments
$\gamma_1$ and $\gamma_2$ 
with $\gamma_1 (0)=p=\gamma_2(0)$, the weak angle $\angle_p^w(\gamma_1,\gamma_2)$ exists.
Then $X$ is 
a $\operatorname{CAT}(0)$-space.
\end{Pro} 
\begin{proof}
We first show that the weak angle satisfies the four axioms
of an angle as formulated in
\cite{brh} II.1.8 (p. 162). Thus we have to show

(A1) $\angle_p^w(\gamma_1,\gamma_2)=\angle_p^w(\gamma_2,\gamma_1)$

(A2) $\angle_p^w(\gamma_1,\gamma_3)\le \angle_p^w(\gamma_1,\gamma_2)+\angle_p^w(\gamma_2,\gamma_3)$

(A3) if $\ga_2$ is the restriction of $\ga_1$ to an initial segment, then 
$\angle_p^w(\gamma_1,\gamma_2)=0$

(A4) if the concatenation of $\ga_1=[p,x]$ and $\ga_2=[p,y]$ is a geodesic $[x,y]$,
then $\angle_p^w(\gamma_1,\gamma_2)=\pi$

Now (A1), (A3), (A4) are trivially true by the definition.
Since weak angles between rays 
coincide (by Lemma \ref{lemma-rays-ultrarays-angles}) with the angles of 
their associated ultrarays and since such angles
satisfy the triangle inequality also (A2) holds.
Furthermore we see that $\angle_p^w([p,x],[p,y]) \le \angle_p(x,y)$,
which follows from the Busemann convexity.
Now implication (4) $\Longrightarrow$ (2) on p.161 of II.1.7 in \cite{brh} remains
valid, if one replaces the Alexandrov angle in condition (4) through the weak angle
(cf. the paragraph after II.1.8 in \cite{brh}).
This implies that $X$ is a $\operatorname{CAT}(0)$-space.
\end{proof}


\subsection{A convex hull proposition}

The purpose of this subsection is to prove 
\begin{Pro} \label{prop-convex-hull} (compare \cite{r}, p. 432, par. 7 and p. 463, par.20 as well as
\cite{bow}, Lemma 1.1. and the remark after its proof)
Let $X$ be Busemann convex and let $\gamma_1,\gamma_2:I\longrightarrow X$ be two
linearly reparameterized (finite or infinite) geodesics in $X$ such that 
$t\mapsto |\gamma_1(t)\gamma_2(t)|$ is affine. Then the convex hull of $\gamma_1$ and $\gamma_2$
is a convex subset of a two-dimensional normed vector space.
\end{Pro}

Given a geodesic metric space $X$, a function $f:X\longrightarrow \mathbb{R}$ is called {\it affine}
if its restriction to each affinely parameterized geodesic $\gamma$ in $X$ satisfies
$f(\gamma (t))=at+b$ for some numbers $a,b\in \mathbb{R}$ that may depend on $\gamma$.
We say that affine functions on $X$  {\it separate points}, if for each pair of distinct points
$x,x'\in X$ there is an affine function $f:X\longrightarrow \mathbb{R}$ with $f(x)\neq f(x')$.
With this terminology the following theorem has been proven in \cite{hl}.
\begin{Thm} \label{theo-affine-separate} (Theorem 1.1 in \cite{hl})
Let $X$ be a geodesic metric space. If affine functions on $X$ separate points then
$X$ is isometric to a convex subset of a normed vector space with a strictly convex norm.
\end{Thm}
Using this result, we are able to provide the \\

{\bf Proof of Proposition \ref{prop-convex-hull}:}
Let $y_t:I_t\longrightarrow X$ be the geodesic from $\gamma_1(t)$ to $\gamma_2(t)$
where $I_t=[0,|\ga_1(t)\ga_2(t)|]$ .
Let $C_0:=\bigcup_{t\in I} y_t(I_t)$. 

{\it Claim:} $C_0$ is convex, i.e., $C_0=C$. \\
In order to prove this claim, we may assume that $\gamma_1$ and $\gamma_2$ are closed. Then $C_0$ is closed
and it is sufficient to prove that for $t_1,t_2\in I$, $s_1\in I_{t_1}$ and $s_2\in I_{t_2}$ the midpoint 
of $y_{t_1}(s_1)$ and $y_{t_2}(s_2)$ is contained in $C_0$. \\
Let $m$ denote the unique midpoint map on $X$ and set
\begin{displaymath}
m_0 \, := \, m(y_{t_1}(s_1),y_{t_2}(s_2)), \;\;\; 
m_1 \, := \, \gamma_1(\frac{t_1+t_2}{2}), \;\;\; \mbox{and} \;\;\; m_2 \, := \, \gamma_2(\frac{t_1+t_2}{2}).
\end{displaymath}

\begin{figure}[htbp]
\psfrag{m0}{$m_0$}
\psfrag{m1}{$m_1$}
\psfrag{m2}{$m_2$}
\psfrag{m3}{$m_3$}
\psfrag{m4}{$m_4$}
\psfrag{g1}{$\gamma_1$}
\psfrag{g2}{$\gamma_2$}
\psfrag{g1t1}{$\gamma_1(t_1)$}
\psfrag{g1t2}{$\gamma_1(t_2)$}
\psfrag{g2t1}{$\gamma_2(t_1)$}
\psfrag{g2t2}{$\gamma_2(t_2)$}
\psfrag{yt1s1}{$y_{t_1}(s_1)$}
\psfrag{yt2s2}{$y_{t_2}(s_2)$}
\includegraphics[width=0.9\columnwidth]{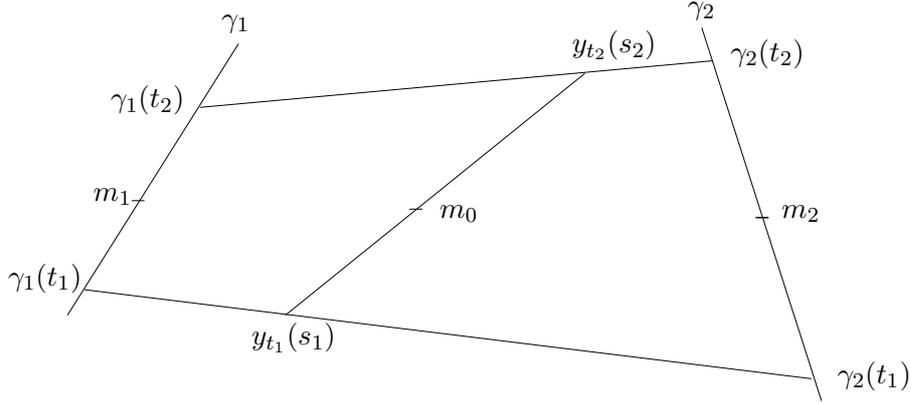}
\caption{This figure visualizes our notation used in the proof of Proposition \ref{prop-convex-hull}.}
\end{figure}

From the convexity of the distance function and the fact that
\begin{displaymath}
|m_1m_2| \; = \; \frac{1}{2} \Big[ |\gamma_1(t_1)\gamma_2(t_1)| \; + \; |\gamma_1(t_2)\gamma_2(t_2)| \Big] ,
\end{displaymath}
we deduce that $|m_1m_2|=|m_1m_0|+|m_0m_2|$.
Thus $m_0$ is contained in the geodesic $y_{\frac{t_1+t_2}{2}}$. This proves the claim. \\

In fact, with the same reasoning as above we deduce more, namely that $|m_0m_1|=\frac{s_1+s_2}{2}$.
This shows that the function $F_1:C\longrightarrow \mathbb{R}$, given through
$F_1(y_t(s))=s$ is affine on $C$. \\
Moreover, the function $F_2:C\longrightarrow \mathbb{R}$, given through $F_2(y_t(s))=t$ is affine as well, 
since the midpoint between $y_{t_1}(s_1)$ and $y_{t_2}(s_2)$ lies on the geodesic $y_{\frac{t_1+t_2}{2}}$. \\
Now the affine functions $F_1$ and $F_2$ separate the points of $C$. Thus we can apply 
Theorem \ref{theo-affine-separate}, which finishes the proof. 
\hfill$\Box$


\subsection{Proof of Theorem \ref{theo-cat0-busemann-ptolemy}}

Every $\operatorname{CAT}(0)$-space is both, Ptolemy and Busemann convex. It remains
to show that a Busemann convex and Ptolemy metric space is already $\operatorname{CAT}(0)$. \\
In order to reach a contradiction, suppose that $X$ is Ptolemy and Busemann convex but not
$\operatorname{CAT}(0)$. Then, due to Proposition \ref{prop-cat0-angles} 
there do exist two geodesics $\gamma_1$ and $\gamma_2$ that do not enclose a weak angle
at their common starting point $p=\gamma_1(0)=\gamma_2(0)$. Let $\overline{\gamma_1}$ and 
$\overline{\gamma_2}$
denote the geodesic rays defined in $Y=\lim_\omega \{ nX,p\}_n$ as above. Then, due to 
Lemma \ref{lemma-rays-ultrarays-angles},
$\overline{\gamma_1}$ and $\overline{\gamma_2}$ do not enclose an angle in $\{ p\}_n\in Y$ either. \\
Now $Y$ is Busemann convex 
by Corollary \ref{cor-often-ptolemy:busemann}. Moreover, the function
$t\mapsto |\overline{\gamma_1}(t) \overline{\gamma_2}(t)|$ is linear. Thus, by Proposition \ref{prop-convex-hull},
the convex hull $C$ of $\overline{\gamma_1}$ and $\overline{\gamma_2}$ is isometric
to a convex set of a two-dimensional normed vector space. 
Since $Y$ is Ptolemy, $C$ is 
flat by Corollary \ref{cor:schoenberg}.
It follows
that $\overline{\gamma_1}$ and $\overline{\gamma_2}$ enclose an angle, which yields the desired contradiction.
\hfill $\Box$



\end{document}